\input amstex
\documentstyle{amsppt}
\magnification = \magstephalf
\NoBlackBoxes
\NoRunningHeads
\TagsOnRight
\def\metric{$\langle\cdot\, ,\cdot\rangle$}
\def\CC{{\Bbb C}}
\def\RR{{\Bbb R}}
\def\Tr{\operatorname{Trace}}
\def\oTr{\operatorname{\overline{Trace}}}
\def\Trn#1{\Tr_N\left({#1}_{(s)}\right)}
\def\End{\operatorname{End}}
\def\Im{\operatorname{Im}}
\def\Ker{\operatorname{Ker}}

\def\res{\operatorname{res}}
\def\Res{\operatorname{Res}}
\def\oRes{\operatorname{\overline{Res}}}
\def\Id{\operatorname{Id}}
\def\calp{\Cal P}
\def\ocalp{\overline{\Cal P}}

\topmatter
\title 
An Extension of the Work of V. Guillemin on Complex Powers and Zeta
Functions of Elliptic Pseudodifferential Operators
\endtitle
\author
Bogdan Bucicovschi
\endauthor
\affil
Ohio State University
\endaffil
\address
Department of Mathematics,
231 W 18th Ave.,
Columbus, OH 43210
\endaddress
\email
bogdanb\@math.ohio-state.edu
\endemail
\abstract
The purpose of this note is to extend the results of Guillemin in [G] on
elliptic self-adjoint pseudodifferential operators of order one, from
operators defined on smooth functions on a closed manifold (scalar 
operators) to operators defined on smooth sections in a vector bundle
of Hilbert modules of finite type over  a finite von~Neumann algebra. 
\endabstract
\endtopmatter
\document
\heading
0. Introduction
\endheading

Let $M$ be a closed Riemannian manifold of dimension $m$ and $E$ a
 vector bundle over $M$ endowed with a hermitian metric. The fibers of
$E$ are finite dimensional vector spaces over $\CC$ or, more general,
finite type Hilbert modules over a von~Neumann algebra $\Cal A$. The
 first situation corresponds to the case $\Cal A=\CC$.
Throughout this paper we will denote by $\Psi (E)$ or simply by $\Psi$ the 
algebra of 
classical pseudodifferential operators acting on smooth sections in
 $E$ (for the case when $\Cal A$ is an arbitrary von~Neumann algebra,
 see [BFKM] for definitions and properties). We 
will also denote by $\Psi^s (E)$ the subspace of pseudodifferential operators 
of complex order $s$. The total symbol $\sigma_{\text{total}}(x,\xi)$ of such an 
operator $A \in \Psi^s$ has locally an asymptotic expansion of the form:
$$
\sigma_{\text{total}}(x,\xi) \sim \sum_{k=0}^\infty\sigma_{s-k}(x,\xi)
$$
where $\sigma_{s-k}(x,\xi)$ are sections of the endomorphism bundle of the pull-back 
of $E$ with respect to the projection map $T^*(M) \setminus \{0\}\to M$.
Each section $\sigma_{s-k}(x,\xi)$ is a homogeneous function in the 
variable $\xi$ of degree of homogeneity $s \in \CC$, $\sigma_{s-k}(x,
\lambda\xi)=\lambda^{s-k}(x,\xi)$ for any $\lambda>0$.
\par
The space $C^{\infty}(E)$ of smooth sections of $E$ over $M$ has a canonical 
metric
$$
\langle f,g\rangle=\int\limits_M\langle f(x),g(x)\rangle_x\, d\,vol
$$
where $\langle\cdot\, ,\cdot\rangle_x$ is the hermitian metric in the fibre above $x\in M$.
The $L^2$ completion of $C^\infty (E)$ with respect to \metric \ will 
be denoted by $L^2 (E)$. A pseudodifferential operator becomes an unbounded 
operator on $L^2(E)$.
\par
We will consider now an elliptic  pseudodifferential operator of order one 
$A \in \Psi ^1$ which is self-adjoint and positive with respect to \metric. 
Suppose that the spectrum of $A$ is included in the interval
$(\epsilon, \infty)$ for a sufficiently small $\epsilon >0$. Then one can 
define the complex powers $A^s$, $s\in \CC$ in the following way
$$
A^s = \frac 1 {2\pi i} \int_\gamma\lambda^s(\lambda-A)^{-1}\, d\lambda \qquad 
\text{when $Re(s)<0$} \tag 0.1
$$
(where $\gamma$ is a contour in the complex plane obtained by joining two lines
parallel to the negative real axis by a circle around the origin) and
$$
\qquad A^s = A^{s-k} A^k \qquad \text{for $Re(s)\geq 0$} \tag 0.2
$$
for large enough $k \in \Bbb Z$ so that $s-k$ is negative.
\par
One of the goals of our paper is to show that $A^s$ is a 
pseudodifferential operator of complex order $s$. We remind the reader
that this fact has been proven first by Seeley [S] in the case of finite 
dimensional hermitian bundle $E$ and
extended to the case of von~Neumann bundles in [BFKM]. We will follow a
different approach due to Guillemin [G].
In the same spirit of Guillemin, we will show that the zeta function of $A$ defined as:
$$
\zeta_A(s) = \Tr_N (A^s) \qquad \text{for $Re(s) < -m$}
$$
has a meromorphic extension over the complex plane $\CC$ with at most 
simple poles at $-m$, $-m+1$, \dots . The residue of $\zeta_A$ at $-m$ will 
be equal to a quantity that depens only on the principal symbol $\sigma_1$ of 
the operator $A$.
\par
Guillemin treatment in [G] covers the case of 
pseudodifferential operators acting on smooth functions on $M$. We will
extend his methods to the case of sections in the vector bundle $E$. The 
main difficulty arises from the fact that the algebra of endomorphisms of 
$E$ is noncommutative (fiberwise it is equal to the algebra of the
$\Cal A-$invariant endomorphisms of the fiber, as compared to
Guillemin's case where the fiber is canonically $\Bbb C$ ).
Our paper has two main sections. In the first part we will show 
that $A^s$ is a 
pseudodifferential operator of order $s \in \CC$. The second part will be 
devoted to the zeta function of $A$.
\par
Throughout the paper $A$ will be a classical pseudodifferential operator of 
order~$1$. The case of an operator of any other positive order can
be reduced to the case in which the order is equal to $1$. 
\par
I thank prof.~D.~Burghelea for his support and suggestions  and
prof.~T.~Kappeler and L.~Friedlander for useful discussions.

\vskip 1 cm

\heading
1. Complex Powers of Pseudodifferential Operators
\endheading

The goal of the section is proving the following:
\proclaim 
{Theorem 1.1}
Let A be a positive, self-adjoint pseudodifferential operator of order
one. Suppose that $Spec(A)\in (\epsilon,\infty)$ for a sufficiently
small $\epsilon >0$. Then its complex powers $A^s$, defined as in
(0.1) and (0.2), are pseudodifferential operators of order $s \in \CC$.
\endproclaim
To show this we will need the following:
\proclaim
{Proposition 1.2}
There exists a holomorphic family of pseudodifferential operators $A_s$
for $s \in \CC$ such that $A_0 = Id$, $A_sA_t = A_{s+t}$ and the difference
$A_1 - A$ is a smoothing operator.
\endproclaim

$(A_s)_{s \in \CC}$ can be thought of as an approximation of the powers 
of $A$ that lie inside $\Psi$. We will show that $A_s - A^s$ are
smoothing operators. Then Theorem 1.1 becomes a straightforward
corollary of Proposition 1.2.
\par
To construct the family $(A_s)_{s \in \CC}$ it will be convenient to
consider the cohomology of 
the group $(\CC,+)$ with coefficients in the representation of
$(\CC,+)$ on the space of sections $C^\infty(\End(\tilde E))$. Here 
$\tilde E$ is the pull-back of the initial vector bundle $E$ over M
with respect to the projection map of the cosphere bundle $ S^*(M) \to
M$. This construction generalizes the cohomology
considered by Guillemin in [G] for the trivial representation of
$(\CC,+)$ on the space of smooth functions on $S^*(M)$.
\par
Let $\sigma$ be a fixed section $\sigma :S^*(M) \to \End(\tilde E)$ so that $\sigma(x,\xi):
E_x \to E_x$ is an invertible positive self-adjoint endomorphism for 
any $(x,\xi) \in S^*(M)$ ($\sigma$ will  be the restriction of the
principal symbol of $A$ to $S^*(M)$). The representation of $(\CC,+)$
on $C^\infty(\End(\tilde E))$ we consider is the following one:
any $s \in \CC$ acts on a section $g:S^*(M) \to \End(\tilde E)$ by $s 
\cdot g =\sigma^{-s}g\sigma^s$. 
\par
Let ${\Cal C}^r={\Cal C}^r(\CC; C^\infty(\End(\tilde E))$ be the space of functions 
$$
f:\undersetbrace \text{$r$ times} \to {\CC\times\CC\times\dots\CC} \to 
C^\infty(\End(\tilde E))
$$ 
that are smooth, $f(\cdot)(x,\xi):\CC\times\CC\times\dots\CC \to \End(E_x)$ 
are holomorphic for any fixed $(x,\xi) \in S^*(M)$ and $f(s_1,\dots,s_r)=0$ if 
at least one $s_i$ is equal to zero.
 
Let $\delta^r:{\Cal C}^r \to {\Cal C}^{r+1}$ defined as:
$$
\multline
(\delta ^r f)(s_0,s_1,\dots,s_r) = s_0\cdot f(s_1,\dots,s_r)+\sum_{i=1}^r (-1)^i
f(s_0,\dots,s_{i-1}+s_i,\dots,s_r)\\ +(-1)^{r+1}f(s_0,\dots,s_{r-1}).
\endmultline
$$
Let ${\Cal H}^r(\CC ; C^\infty(\End(\tilde E))={\Ker \delta^r}/{\Im \delta^{r-1}}$.
\proclaim{Proposition 1.3}
${\Cal H}^2(\CC;C^\infty(\End(\tilde E))=0$.

Moreover, for each 2-cocycle $f$ there exists a unique 1-cochain $h$ such 
that $\delta h = f$ and $h$ has a prescribed value at $1$, $h(1)$.
\endproclaim
\demo{Proof}
Let $f:\CC\times\CC \to C^\infty(\End(\tilde E))$ so that for all $a,b,c \in \CC$ 
$$
\left\{\aligned
&f(0,b)=f(a,0)=0 \\
&(\delta^2 f)(a,b,c) = a \cdot f(b,c)-f(a+b,c)+f(a,b+c)-f(a,b)=0
\endaligned
\right.
$$ 
We will try to find $h:\CC \to C^\infty(\End(\tilde E))$ such
that 
$$
(\delta^1 h)(a,b)= \sigma^{-a} h(b) \sigma^a -h(a+b)+h(a) = f(a,b)
$$

The existence of an $h$ as above implies:
$$
h'  (a) = \sigma^{-a} h'  (0) \sigma^a - \frac {\partial f}{\partial b}(a,0)
\tag 1.1
$$
Consider $h$ to be the unique solution of the previous equation with $h(0)=0$ 
and with a fixed prescribed value at $1$, $h(1)$. $h$ can be found
in the following way:\newline
Let $\Phi(t)$ be the automorphism of $C^{\infty}(\End(\tilde E))$
given by $A\to\sigma^{-t}A\sigma^t$. Then
$$
h(a)=-\int_0^a\frac {\partial f}{\partial b}(t,0)\,dt +\int_0^a
\Phi(t)(h'(0))\,dt
$$
If $\displaystyle T(a)A=\int_0^a \Phi(t)A\,dt$, then, in order to
get any prescribed value for $h(1)$, we need to show that $T(1)$ is
surjective. Indeed, we have:
$$
\align
T(1)A&=\int_0^{\frac 1 2}\sigma^{-t}A\sigma^t \,dt +\int_{\frac 1 2}^1
\sigma^{-t}A\sigma^t \,dt \\
&=T({\tfrac 1 2})A +\Phi({\tfrac 1 2})T({\tfrac 1 2})A=(\Id+\Phi({\tfrac1 2}))
T({\tfrac 1 2})A\\
\intertext{and by induction}
T(1)A&=(\Id+\Phi({\tfrac 1 2}))(\Id+\Phi({\tfrac 1 4}))\dots(\Id +\Phi
({\tfrac{1}{2^n}}))T({\tfrac 1{2^n}})A
\endalign
$$
But the map $\displaystyle A\to 2^n\!\int_0^{\frac 1{2^n}}\sigma^{-t}A\sigma^t\, dt$ is close to
the identity for a sufficiently large $n$ so $T({\frac 1{2^n}})$ is
invertible. $(\Id +\Phi({\tfrac 1{2^i}}))$ is invertible as well, 
because $\Phi(t)$ is positive self-adjoint for any real $t$.

Thus we obtain a continuous map $h:\CC \to C^\infty(\End(\tilde E))$ that is 
holomorphic in all fibers $E_{(x,\xi)}$, $h\in \Cal C^1$. We will show that $\delta \, h = f$ 
so $f$ is a coboundary. To see this, let 
$$
g(a,b) = f(a,b)-(\sigma^{-a} h(b) \sigma^a - h(a+b) + h(a))
$$
Clearly $\delta \, h =f$ if and only if $g \equiv 0$. Denote by $\dfrac 
{\partial}{\partial b}$ the partial derivative with respect to the
second variable. Then:
$$
\frac {\partial g}{\partial b}(a,b)=\frac {\partial f}{\partial b}
(a,b)-\sigma^{-a}h' (b)\sigma^a+h'(a+b) \tag 1.2
$$
From (1.1) we get:
$$
\align
h'(b)&=\sigma^{-b}h(0)\sigma^b - \frac {\partial f}{\partial b}(b,0) \quad 
\text{and}\\
h'(a+b)&=\sigma^{-(a+b)}h'(0)\sigma^{(a+b)}-\frac {\partial f}{\partial b}
(a+b,0)
\endalign
$$
These two equalities and (1.2) imply
$$
\align
\frac {\partial g}{\partial b}(a,b)&=\frac {\partial f}{\partial b}(a,b)- 
\sigma^{-a}\left(\sigma^{-b}h'(0)\sigma^b - \frac {\partial f}{\partial b}(b,0)
\right)
\sigma^a+\sigma^{-(a+b)}h'(0)\sigma^{(a+b)}-\\
&\quad -\frac {\partial f}{\partial b}
(a+b,0)= \\
&= \sigma^{-a}\frac {\partial f}{\partial b}(b,0)\sigma^a-\frac {\partial f}
{\partial b}(a+b,0)+\frac {\partial f}{\partial b}(a,b)=\\
&=\frac {\partial }{\partial c}\left[(\delta^2 \, f)(a,b,c)\right]_{|_{c=0}}
\endalign
$$

So $\dfrac {\partial g}{\partial b} = 0$ hence $g(a,b)$ is constant in $b$. 
When $b=0$ we have 
$$
g(a,0)=f(a,0)-\left(\sigma^{-a}h(0)\sigma^a-h(a)+h(a)\right)=0
$$
So $g\equiv 0$. Because $f$ was chosen arbitrarily we conclude 
${\Cal H}^2(\CC ; C^\infty(\End(\tilde E))=0$.

\qed
\enddemo

We now proceed with the proof of Proposition 1.2. To show the existence of a family
$(A_s)_{s\in \CC}$ as stated in the proposition, we will show that there exists a 
family $(A_{(s)})_{s\in \CC}$ of pseudodifferential operators that satisfies 
the conditions of Proposition 1.2  only up to smoothing operators. 
More precisely:
\proclaim{Proposition 1.4}
There exists a holomorphic family of pseudodifferential operators $(A_{(s)})_
{s \in \CC}$ with principal symbols $\sigma_{pr}(A_{(s)})=(\sigma_{pr}(A))^s$ 
such that $A_{(0)}=Id$, $A_{(1)} \equiv A$ and $A_{(s)}A_{(t)} 
\equiv A_{(s+t)}$ modulo smoothing operators. This family is unique up to 
smoothing operators.  
\endproclaim
\demo{Proof}
The statement of the Theorem is equivalent to:
$$
\cases A_{(s)}A_{(t)}A^{-1}_{(s+t)} \equiv Id  &\pmod {\Psi^{-\infty}} \\
A_{(1)}A^{-1} \equiv Id &\pmod {\Psi^{-\infty}} \\
A_{(0)}= Id&
\endcases
\tag 1.3
$$
(we denoted the space of smoothing operators by $\Psi^{-\infty}$)\newline
To prove Proposition 1.4, we will construct $A_{(s)}$ inductively in
$k \in \Bbb N$, such that
$$
\cases A_{(s)}A_{(t)}A^{-1}_{(s+t)} \equiv Id  &\pmod {\Psi^{-k}} \\
A_{(1)}A^{-1} \equiv Id &\pmod {\Psi^{-k}} \\
A_{(0)}= Id&
\endcases
\tag 1.4
$$

For $k=1$ we can choose $(A_{(s)})_{s\in \CC}$ to be a holomorphic family of 
pseudodifferential operators of order $s$ with the principal symbol equal to
$\sigma^s$ where $\sigma$ is the principal symbol of $A$. We can construct 
such a family using a partition of unity. Moreover $A_{(0)}$ can be chosen to 
be the identity. The operators $A_{(s)}A_{(t)}A^{-1}_{(s+t)}$ and $A_{(1)}
A^{-1}$ are operators of order $0$ with the principal symbol equal to the 
principal symbol of the identity. The relations (1.4) are satisfied modulo 
$\Psi^{-1}$.

Now suppose that the relations (1.4) hold for a certain $k \in \Bbb N$. We will
construct a new family $(\tilde A_{(s)})_{s \in \CC}$ that satisfies (1.4) for
$k+1$, that is of the following form:
$$
\tilde A_{(s)} = A_{(s)}(Id - H_{(s)}), \qquad H_{(s)} \in \Psi^{-k}\tag 1.5
$$
In this way $\tilde A_{(s)}-A_{(s)} \in \Psi^{s-k}$. We have:
$$
\align \tilde A_{(s)} \tilde A_{(t)}\tilde A_{(s+t)}^{-1} &\equiv 
A_{(s)}(Id - H_{(s)})A_{(t)}(Id - H_{(t)})(Id + H_{(s+t)})A_{(s+t)}^{-1}
\equiv \\
&\equiv A_{(s)}A_{(t)}A_{(s+t)}^{-1} - A_{(s)}H_{(s)}A_{(t)}A_{(s+t)}^{-1} -
A_{(s)}A_{(t)}H_{(t)}A_{(s+t)}^{-1} +\\
&\quad + A_{(s)}A_{(t)}H_{(s+t)}A_{(s+t)}^{-1} \\
&\equiv Id +F_{(s,t)}-A_{(s)}H_{(s)}A_{(t)}A_{(s+t)}^{-1}-A_{(s)}A_{(t)}H_{(t)}
A_{(s+t)}^{-1}+\\
&\quad + A_{(s)}A_{(t)}H_{(s+t)}A_{(s+t)}^{-1}
\qquad \pmod {\Psi ^{-k-1}} \tag 1.6
\endalign
$$ 
where $F_{(s,t)}=A_{(s)}A_{(t)}A_{(s+t)}^{-1}-Id$ , $F_{(s,t)} \in 
\Psi^{-k}$ by the induction step. To proceed with the induction we have to 
find a family $(H_{(s)})_{s \in \CC}$ that makes the right hand~side of the 
equivalence (1.6) equal to the identity modulo $\Psi^{-k-1}$. If $\sigma_{pr}
(F(s,t))$ and $h(s)=\sigma_{pr}(H(s))$ are the principal symbols, then the 
condition on $H(s)$ is equivalent to:
$$
\gather
\sigma_{pr}(F(s,t))=\sigma^sh(s)\sigma^{-s} + \sigma^{s+t}h(t)\sigma^{-(s+t)}
-\sigma^{s+t}h(s+t)\sigma^{-(s+t)}  \quad \text{or}\\
\sigma^{-(s+t)}\sigma_{pr}(F(s,t))\sigma^{s+t} = \sigma^{-t}h(s)\sigma^t -
h(s+t) +h(t) \tag 1.7
\endgather
$$
Because both sides are sections in the bundle $\End(\tilde E)$ over $T^*(M)\setminus\{0\}$
of degree of homogeneity $-k$, then the above equality is satisfied if it holds
when both sections are restricted to the cosphere bundle $S^*(M)$. Let:
$$
f(t,s)=\sigma^{-(s+t)}\sigma_{pr}(F(s,t))\sigma^{s+t} \qquad \text{restricted 
to }S^*(M) \tag 1.8
$$
We will show that $f \in {\Cal C}^2(\CC; C^\infty(\End(\tilde E))$ and $\delta^2f=0$.
Then $h$ as in (1.7) will be a 1-cochain so that $\delta h = f$. 
\par
We would also want the second condition of (1.4) to be satisfied so:
$$
\align
A^{-1}\tilde A_{(1)} &\equiv A^{-1}A_{(1)}(Id - H_{(1)}) \equiv \\
&\equiv Id +(A^{-1}A_{(1)}- Id)-A^{-1}A_{(1)}H_{(1)} \\
&\equiv Id \qquad \pmod{\Psi^{-k-1}}
\endalign
$$ 
and this holds if 
$$
h(1)=\sigma_{pr}(A^{-1}A_{(1)}-Id)\tag 1.9
$$
(we already know that $(A^{-1}A_{(1)} - Id) \in \Psi^{-k}$ from the induction 
step).
\par
We will have to show that $f$ is a cocycle in $\Cal C^2$. Obviously, $f(0,t)=
f(s,0)=0$. We have:
$$
\align
(\delta^2f)(s,t,r)&=\sigma^{-s}f(t,r)\sigma^s - f(s+t,r) + f(s,t+r) -f(s,t)=\\
&=\sigma^{-s}\left[\sigma^{-(t+r)}\sigma_{pr}(F(r,t))\sigma^
{t+r}\right]\sigma^s - \sigma^{-(s+t+r)}\sigma_{pr}(F(r,s+t))\sigma^{s+t+r} \\
&\quad+\sigma^{-(s+t+r)}\sigma_{pr}(F(t+r,s))\sigma^{s+t+r}-
\sigma^{-(s+t)}\sigma_{pr}(F(t,s))\sigma^{s+t}=0
\endalign
$$
is equivalent to
$$
\sigma_{pr}(F(r,t)) - \sigma_{pr}(F(r,s+t)) + \sigma_{pr}(F(t+r,s)) - 
\sigma^r\sigma_{pr}(F(t,s))\sigma^{-r}=0 \tag 1.10
$$
To see this, consider the following equivalences modulo $\Psi^{-k}$:
$$
\align
&\quad(Id+F(r,t))(Id+F(t+r,s))(Id-F(r,s+t))A_{(r)}(Id-F(t,s))A_{(r)}^{-1} 
\equiv \\
&\equiv A_{(r)}A_{(t)}A_{(t+r)}^{-1}A_{(t+r)}A_{(s)}A_{(s+t+r)}^{-1}
A_{(s+t+r)}A_{(s+t)}^{-1}A_{(r)}^{-1}A_{(r)}A_{(s+t)}A_{(s)}^{-1}A_{(t)}
^{-1}A_{(r)}^{-1} \equiv \\
& \equiv Id
\endalign
$$
and the first term is also equivalent to
$$
Id+F(r,t)-F(r,s+t)+F(t+r,s)-A_{(r)}F(t,s)A_{(r)}^{-1}
$$
which proves (1.10). So $f(s,t)=\sigma^{-(s+t)}\sigma_{pr}(F(t,s))\sigma^{s+t}$
is a cocycle. 
\par
Proposition 1.3 provides us with a family $h(s)$ such that 
$\delta h=f$. We can choose this family so that (1.9) holds as well. This 
determines $h$ in a unique way. If $(H_{(s)})_{s \in \CC}$ is a holomorphic 
family of pseudodifferential operators of fixed order $-k$ with principal 
symbol $h(s)$ and $H_{(1)}=Id$, then $\tilde A_{(s)} = A_{(s)}(Id-H_{(s)})$ 
satisfies the equivalences (1.4) modulo $\Psi^{-k-1}$.
\par
In this way we obtain a sequence of families of operators $(A_{(s)}^{(k)})_{s 
\in \CC}$ that satisfy the relations (1.4) for each $k \in \Bbb N$. Moreover, 
$A_{(s)}^{(k+1)}-A_{(s)}^{(k)} \in \Psi ^{s-k}$. Then, using a
standard procedure as in
Lemma 1.2.8 in [Gi], we can construct a family $(A_{(s)})_{s \in \CC}$ whose
asymptotic expansion of the total symbol will be equal to:
$$
\sigma_{\text{total}}(A_{(s)})\sim \sigma_{\text{total}}(A_{(s)}^{(1)})+\sum_{k \ge 0}\sigma_{\text{total}}(A_{(s)}^{(k+1)}-A_{(s)}^{(k)}) 
$$ 
The family $(A_{(s)})_{s \in \CC}$ will satisfy the conditions of
Proposition 1.4. 
\par
$(A_{(s)})_{s \in \CC}$ is unique up to smoothing operators because it must 
satisfy the relations (1.4) for all $k \in \Bbb N$ and so it must be equal 
to $(A_{(s)}^{(k)})_{s \in \CC}$ modulo $\Psi^{-k}$.

\qed
\enddemo

\demo{Proof of Theorem 1.1 and Proposition 1.2}
Once we obtained the family of pseudodifferential operators
$(A_{(s)})_{s \in \CC}$, the proofs of Thm. 1.1 and Prop. 1.2 are identical to
the proof of Theorem~5.1 in [G]. We can construct the one parameter
group of operators as in Prop 1.2 using the differential equation:
$$
\dot A_s=PA_s \qquad \text{with}\quad A_0=Id
$$ 
where $P=\dot A_{(0)}$. If $A_{(s)}$ is made a selfadjoint family in
$s$ (i.e. $A_{(s)}^*=A_{(\overline{s})}$) by replacing it with 
$\frac 1 2 (A_{(s)}+A_{(\overline{s})}^*)$, $P$ becomes 
a selfadjoint operator. By construction $A_s \in \Psi^s$. 
Then, using a theorem of Stone (Thm VIII.7 and Thm VIII.8 [RS]), it can
be shown that $A_s=(A_1)^s$ with $P$ the infinitesimal generator of
this one parameter group.
In this case:
$$
(A_1)^{\,s}-A^s = \frac 1{2\pi i}\int_{\gamma}\lambda^s(\lambda-A)^{-1}
(A-A_1)(\lambda-A_1)^{-1}\,d\lambda
$$
and this is a smoothing operator. Because $(A_1)^s=A_s \in\Psi^s$ we
obtain $A^s\in\Psi^s$.

\qed
\enddemo

\vskip 1 cm

\heading
2. Zeta Function of an Elliptic Pseudodifferential Operator
\endheading

Let $(A_{(s)})_{s \in \CC}$ be a family of pseudodifferential operators 
depending holomorphically on the complex parameter $s$, $A_{(s)} \in
\Psi^s$. For $Re(s) < - dim(M)$, $A_{(s)}$ is a trace-class operator.

\definition{Definition 2.1} The trace function of the
family$A_{(s)}$ is the holomorphic function 
$\Trn A \quad \text{for }Re(s)<-dim(M).$
 
The von~Neumann trace of $A_{(s)}$ is obtained by integrating the von~Neumann
trace of the Schwartz kernel on $M$ for $Re(s)<-dim(M)$.
If $A$ is an elliptic positive self-adjoint pseudodifferential operator of 
order $1$ with $Spec(A)\in(\epsilon,\infty)$ then its zeta function
$\zeta_A$ is equal to the trace function associated with the family 
of its complex powers $A^s$.
\enddefinition

In this section of our paper we will show that $\Trn A$ has a
meromorphic continuation to the whole complex plane with at most
simple poles at $-m$, $-m+1$, \dots, where $m=dim(M)$. This fact has
been proved by Seeley [S]. Guillemin has a different proof in [G] that
applies only for scalar pseudodifferential operators. We will adapt
his proof for the case of operators that act on sections in a vector 
bundle $E$ over the base space $M$.
\par
We start by recalling some definitions and constructions in [G].
\par
Let $\omega$ be the canonical symplectic form on the cotangent space 
$Y=T^*(M)\setminus\{0\}$. The multiplicative group $(\RR^+,\cdot)$
acts on $Y$ by multiplication along the fibre $(t,(x,\xi))\overset
\rho \to{\to} (x,t\xi)$. By identifying the groups $(\RR^+,\cdot)$
and $(\RR , +)$ via $\ln:\RR^+ \to \RR$, $\rho$ can 
be seen as a 1-parameter group of
isomorphisms. Let $\Xi$ be the vector field on $Y$ associated with
this 1-parameter group and $\alpha=\iota_{\Xi}\,\omega$ be the
contraction of $\omega$ along $\Xi$. Then the $(2m-1)$-form on $Y$, 
$\mu=\alpha\wedge\omega^{m-1}$, is homogeneous of degree $m$, 
$\rho^*_t\,\mu=t^m\mu$, and it is horizontal with respect to the
fibration\linebreak $Y=T^*(M)\setminus\{0\}\overset \pi \to{\to}S^* (M)$.
\par

Let $\Cal B$ be a von~Neumann algebra. In our case, $\Cal B$ will be the
field of complex numbers $\Bbb C$, our initial von~Neumann algebra
$\Cal A$ or $\End_{\Cal A}(V)$, where $V$ is the generic fiber of 
the vector bundle $E\to M$. Let  $\ocalp_s$ be the space of smooth 
homogeneous $\Cal B$-valued functions
defined on $Y$ of degree of homogeneity $s \in \CC$ and $\calp_s$ the
space of smooth scalar functions on $Y$ of degree of homogeneity $s$.
If $f\in\ocalp_{-m}$ then the
$\Cal B$-valued $(2m-1)$ form $f\mu$ is horizontal and
invariant under the action of 
$(\RR^+,\cdot)$ so it is of the form $\pi^*\mu_f$ where $\mu_f$ is a
$(2m-1)$-form on $S^*(M)$.
\definition{Definition 2.2} The residue of $f\in\ocalp_{-m}$ is equal
to the integral
$$
\oRes f=\int_{S^*(M)}\mu_f\quad \in \Cal B
$$
For $f \notin \ocalp_{-m}$ we define $\oRes f=0$.\newline
If $\Cal B = \Bbb C$, we will denote the residue simply by $\Res f$.

\enddefinition

Consider the Poisson bracket $\{\, ,\,\}$ on $T^*(M)$ associated
with the canonical symplectic form $\omega$. Let $\{\calp_s,\ocalp_t
\}$ be the space of functions spanned by $\{f,g\}$ with $f \in
\calp_s$ and $g\in \ocalp_t$. Then $\{\calp_s,\ocalp_t\} \subset
\ocalp _{s+t-1}$. 
Following the same method as in [G] (Theorem 6.2), it can be shown that: 
\par
a) If $s\neq -m$ then $\{\calp_1,\ocalp_s\}=\ocalp_s$.
\par
b) If $s=-m$ then $\{\calp_1,\ocalp_s\}$ 
consists of all functions $f$ for which $\oRes f=0$.
\par
Moreover, one can construct a family of functions $(g_i)_{i\in I}$,
$g_i \in \calp_1$ such that for any analytic family with parameter
$s$, $f_s \in \ocalp_s$, defined on a strip $a-\epsilon \le Im(s)
\le a+\epsilon$, $c\le Re(s) \le d$ for which $\oRes f_{-m}=0$, one
can find $\delta \le \epsilon$ and homogeneous functions $h_{i,s}
\in \ocalp$ which are analytic in $s$ on a narrower strip 
$a-\delta \le Im(s) \le a+\delta$, $c\le Re(s) \le d$, such that 
$$
 \quad f_s=\sum_{i \in I}\{g_i, h_{i,s}\} 
$$
(cf [G], Theorem 6.7)
\par
Let us consider now a holomorphic family of pseudodifferential
operators\linebreak $(A_{(s)})_{s \in \CC}$, $A_{(s)} \in \Psi^s$ and its
associated trace function $\Trn A$. We define the 
residue of the family $A$ to be $\Res A = \Res(\Tr_N\sigma_{pr}
(A_{(-m)}))\in\CC$.
\par
We have the following theorem:
\proclaim{Theorem 2.3} The trace function of the analytic family
$(A_{(s)})_{s\in \CC}$ has a meromorphic continuation to the
whole complex plane with at most simple poles at $-m$,~$-m+1$,~\dots.
The residue of $\Trn A$ at $s=-m$ is equal to
$$
\res_{\mid_{s=-m}}\Trn A =\gamma_0\Res A
$$
where $\gamma_0$ is a constant depending only on $dim(M)$. For 
$A_{(s)}=A^s$ -- the complex powers of an elliptic positive self-adjoint 
pseudodifferential operator of order one, the residue of the 
zeta function at $s=-m$ depends only on its principal symbol 
$\sigma=\sigma_{pr}(A)$ and is equal to $\gamma_0\Res(\sigma^{-m})$.
\endproclaim
\demo{Proof} Let $(U_{\alpha})_{\alpha}$ be an open cover of $M$ with
chosen trivializations of the vector bundle $E$ over each $U_{\alpha}$,
$E_{\mid_{U_{\alpha}}}\cong U_{\alpha}\times V$, with $V$ the generic
fiber. 
Using a partition of unity associated to the open cover 
$(U_{\alpha})_{\alpha}$, we can write:
$$
A_{(s)}=\sum _{\alpha}A_{\alpha\,(s)}+K_{(s)} \tag{2.1}
$$
where $A_{\alpha\,(s)}$ are pseudodifferential operators of order 
$s$ with support inside $U_{\alpha}$ and $K_{(s)}$ is a family of
smoothing operators. Because the residue of the trace function of the
family $A_{(s)}$ and $\Res A$ are both linear in $A$, it is sufficient
to prove the theorem for $A_{\alpha\,(s)}$ and $K_{(s)}$. But 
$K_{(s)}$  is a family of smoothing operators and both the residues of
their trace function and the residue $\Res K$ are zero. Thus we reduced
the proof of the theorem to the case of one family $A_{(s)}=A_{\alpha\,(s)}$ 
supported in an open set $U=U_{\alpha}$ over which we have a 
trivialization of the vector bundle $\chi_{\alpha}:E_{\mid U}\to U\times V$. 
Moreover, because both the trace function $\Trn A$
and the residue $\Res A$ are obtained by integrating quantities that depend
on the local expression of the total symbol of $A_{(s)}$, we can replace 
the bundle $E\to M$ with the trivial bundle $M\times V\to M$ and the 
operators $A_{\alpha\,(s)}$ with the pseudodifferential operators acting
on sections of the trivial bundle $M\times V$ that are supported in
the open set $U_{\alpha}$ and equal to $A_{\alpha\,(s)}$ via the 
isomorphism $\chi_{\alpha}$.
To make things simple, we will denote this new family of operators by 
$A_{(s)}$ as well, and the new trivial bundle by $E$.
\par
Following the ideas in [G], we consider the family $(s+m)A_{(s)}$.
The principal symbol $(s+m)\sigma_{pr}(A_{(s)})$ can be represented by
the $\Cal B$-valued smooth homogeneous functions of degree $s$,
$f_{(s)}:T^*(M)\setminus\{0\}\to \Cal B$, with $\Cal B=\End_{\Cal
A}(V)$. For $s=-m$ we have 
$f=0$, so $\oRes f = 0$. Then there exist $\Cal B$-valued
functions $h_{(s)}^k$, $h_{(s)}^k\in\ocalp_s$ such that
$$
f_{(s)}=\sum_k\{g_k,h_{(s)}^k\}
$$
and $h_{(s)}^k$ are analytic on a strip $a-\epsilon \le Im(s)
\le a+\epsilon$, $c\le Re(s) \le d$.
\par
Let $G_k=G_k'\hat{\otimes} Id$ be a pseudodifferential operator acting
 on the space of sections
 $C^{\infty}(M)\hat{\otimes} V$ of the trivial bundle $E$ with $G_k'$ a scalar 
pseudodifferential operator that  has the
principal symbol equal to $g_k$ and $Id$ the identity operator. Let 
$(H_{(s)}^k)$  be a
holomorphic family of pseudodifferential operators with the principal symbol 
equal to $h_{(s)}^k$. Then the principal symbol of the commutator
is equal to
$$
\sigma_{pr}\left[G_k,H_{(s)}^k\right]=\{g_k,h_{(s)}^k\}
$$
so 
$$
(s+m)A_{(s)}=\sum_k\left[G_k,H_{(s)}^k\right] + B_{(s)}\qquad \text{with }
B_{(s)}\in \Psi^{s-1}.
$$
\par
For $Re(s)$ sufficiently small, $\Tr_N\left([G_k,H^k_{(s)}]\right)=0$, so 
$\Trn A =\break\dfrac 1{s+m}\Trn B$ for $Re(s)<-m$. But $\dfrac 1{s+m}\Trn B$ 
is a meromorphic function on the half-plane $Re(s)<-m+1$ with a simple 
pole at $s=-m$. So $\Trn A$ has a meromorphic extension to 
$Re(s)<-m+1$.
Replacing the family $A_{(s)}$ by $B_{(s)}$ and using an induction argument,
we can extend $\Trn A$ to a meromorphic function on the complex plane 
with at most simple poles at $-m$, $-m+1$,\dots. 
\par
We will compare the residue of $\Trn A$ at $-m$ to the residue of the
family $(A_{(s)})_{s \in \CC}$, $\Res A=\Res (\Tr_N \sigma_{pr}(A_{(-m)}))$. 
Guillemin has showed ([G], Theorem 7.5)
that in the scalar case there exists a constant $\gamma_0$ that depends 
only on the dimension of the manifold $M$ such that 
$$
\res_{\mid_{ s=-m}}\Tr A =\gamma_0\Res A \tag{2.2}
$$ 
We will extend this equality for the pseudodifferential 
operators acting on sections in the vector bundle $E$.
\par
We will show a stronger equality:
$$
\res_{\mid_{ s=-m}}\oTr A=\gamma_0\oRes A_{(-m)} \tag{2.3}
$$ 
where $(A_{(s)})_s$ is a holomorphic family of pseudodifferential
operators acting on the sections of the trivial bundle $M\times V$,
$\oTr A_{(s)}=\int_{M}K_s(x,x)\,dx$ with $K_s(x,y)$ the Schwartz
kernel of $A_{(s)}$, and $\oRes A_{(-m)}=\oRes \sigma_{pr}(A_{(-m)})$,
both sides of the equality $(2.3)$ being in the von~Neumann algebra 
$\Cal B=\End_{\Cal A}(V)$. The equality $(2.2)$ will be then a direct
consequence of $(2.3)$ after passing to the von~Neumann traces.
\par
Both sides of the equality $(2.3)$ depend only on the principal symbol
of the operator $A_{(-m)}$. This is obvious for the right-hand
side. If one considers another family $B_{(s)}$ with
$\sigma_{pr}(B_{(-m)})=\sigma_{pr}(A_{(-m)})$, then $(B_{(s)}-A_{(s)})$ is a
family for which $\oRes \sigma_{pr}(B_{(-m)}-A_{(-m)})=0$, so, by 
a previous observation, $\oTr(B_{(s)}-A_{(s)})$ has a meromorphic
extension which is holomorphic at $s=-m$. So $\oTr B_{(s)}$ and $\oTr A_{(s)}$
will have the same residue at $s=-m$ and this shows that the
left-hand side of $(2.3)$ depends only on $\sigma_{pr}(A_{(-m)})$.
\par
Both sides of $(2.3)$, as functions of holomorphic families, will
factor through the projection ${A_{(s)}}\to \sigma_{pr}(A_{(-m)})
\in\ocalp_{-m}$. It will be sufficient to show that the 
equality $(2.3)$ holds on $\ocalp_{-m}$. 
\par
$\oRes$ vanishes exactly on $\{\calp_1,\ocalp_{-m}\}$ and realizes 
a $\Cal B$ isomorphism\break
$\ocalp_{-m}/\{\calp_1,\ocalp_{-m}\} \overset \sim\,\to\to 
{\Cal B}$. 
For $f\in\{\calp_1,\ocalp_{-m}\}$, $f=\sum \{g_k, h^k\}$, one
can extend it to a holomorphic family of homogeneous symbols of degree
of homogeneity $s\in\CC$ by considering  first the homogenous holomorphic 
extensions $h^k_{(s)}\in\ocalp_s$ and then taking $f_{(s)}=\sum \{g_k,
h^k_{(s)}\}$. If $G_k=G_k'\hat{\otimes} Id$ is a 
pseudodifferential operator such that the scalar operator $G_k'$ has
the principal symbol equal to $g_k$ and $(H^k_{(s)})$ is a holomorphic
family of pseudodifferential operators with the principal symbol equal to
$h^k_{(s)}$, then $A_{(s)}$ defined as $\sum[G_k,H^k_{(s)}]$ has the principal 
symbol at $s=-m$ equal to $f$ and its trace is identically zero. This
shows that $\res_{\mid_{ s=-m}}\oTr A$ vanishes on
$\{\calp_1,\ocalp_{-m}\}$  as well.
Because both $\res_{\mid_{ s=-m}}\oTr A$ and $\oRes A_{(-m)}$ are $\Cal B$
linear, one gets $\res_{\mid_{ s=-m}}\oTr A=\oRes A_{(-m)}\cdot C$ with
$C\in\Cal B$.
\par
Guillemin already showed this equality for a holomorphic family 
of scalar pseudodifferential operators  $(A_{(s)})$ in which
case $C$ is a scalar constant $\gamma_0$. So
$C=\gamma_0\cdot\Id_{\Cal B}$ and the equality $(2.3)$ holds. Passing to
the von Neuman trace, we get $(2.2)$. 

\qed
\enddemo

\enddocument